\documentclass[11pt,bezier]{article}
\pdfoutput=1
\usepackage{amsmath}
\usepackage{amsfonts,amsthm,amssymb}
\usepackage{exscale}
\usepackage{relsize}
\usepackage{graphicx}
\usepackage{epstopdf}
\usepackage{caption}
\usepackage{float}

\usepackage{amssymb}
\usepackage{amsmath}
\newtheorem{lemma}{Lemma}[section]
\newtheorem{theorem}{Theorem}[section]
\newtheorem{corollary}[theorem]{Corollary}

\usepackage{graphicx}
\makeatletter
\newcommand{\Rmnum}[1]{\expandafter\@slowromancap\romannumeral #1@}
\makeatother
\theoremstyle{definition}

\usepackage[justification=centering]{caption}
\allowdisplaybreaks[3]
\begin{document}
\title{The 3-restricted edge-connectivity of the direct product graphs\footnote{The research is supported by National Natural Science Foundation of China (12261086).}}
\author{Wenxin Wang, Yingzhi Tian\footnote{Corresponding author. E-mail: 107552300678@stu.xju.edu.cn (W. Wang); tianyzhxj@163.com (Y. Tian).} \\
{\small College of Mathematics and System Sciences, Xinjiang
University, Urumqi, Xinjiang 830046, PR China}}

\date{}
\maketitle

\begin{sloppypar}

\noindent{\bf Abstract} An edge subset \( S \subseteq E(G) \) is called a 3-restricted edge-cut if \( G - S \) is disconnected and each component of \( G - S \) contains at least three vertices. The 3-restricted edge-connectivity of a graph \( G \), denoted by \( \lambda_3(G) \), is defined as the minimum cardinality among all  3-restricted edge-cuts if there are at least one; otherwise, \( \lambda_3(G) = +\infty \).
It is proved that $\lambda_3(G)\leq\xi_3(G)$ if $G$ has a 3-restricted edge-cut, where $\xi_3(G) = \min \left\{ |[X, V(G) \setminus X]_G|:|X| = 3 \text{ and } G[X] \text{ is connected} \right\}.$
If \( \lambda_3(G) = \xi_3(G) \), then \( G \) is said to be maximally 3-restricted edge-connected. The direct product of two graphs $G$ and $H$, denoted by $G \times H$, is defined as the graph with vertex set \( V(G \times H) = V(G) \times V(H) \), where two vertices \( (u_1, v_1) \) and \( (u_2, v_2) \) are adjacent in \( G \times H \) if and only if \( u_1u_2 \in E(G) \) and \( v_1v_2 \in E(H) \). In this paper, we determine, for a regular connected graph \( G\), the 3-restricted edge-connectivity of \( G \times C_n \), \( G \times K_n \) and \( G \times T_n \), where \( C_n \), \( K_n \) and \( T_n \) are the cycle, the complete graph and the total graph with \( n \) vertices, respectively. As corollaries, we establish sufficient conditions for the direct product graphs \( G \times C_n \), \( G \times K_n \) and \( G \times T_n \) to be maximally 3-restricted edge-connected.


\section{Introduction}

The topological structure of an interconnection network is often modeled by a graph $G = (V(G), E(G))$, where the vertex set $V(G)$ represents the nodes of  the network and the edge set $E(G)$  represents the communication links between nodes. The $order$ of \( G \) is the number of its vertices.  For a vertex \( u \in V(G) \), the $neighborhood$ $N_G(u)$ of \( u \) in \( G \) is defined as \( \{ v \in V(G) \mid v \text{ is adjacent to } u \} \); the $degree$ \( d_G(u) \) of   \( u \) in \( G \)  is the number of edges incident with \( u \), where a loop counts twice. The $minimum$ $degree$ $\delta(G)$ and the $maximum$ $degree$ $\Delta(G)$ of \( G \) are defined as $ \min\{ d_G(u) \mid u \in V(G)\}$ and $ \max\{ d_G(u) \mid u \in V(G) \}$, respectively. If each vertex of $G$ has the same degree $k$, then $G$  is said to be $k$-$regular$.
For a vertex subset $U \subseteq V(G)$, the $induced$ $subgraph$ of $U$ in $G$, denoted by $G[U]$, is the graph with vertex set $U$, where two vertices $u$ and $v$ in $U$ are adjacent in $G[U]$ if and only if they are adjacent in $G$. For an edge \( e = uv \in E(G) \), the $edge$-$degree$ of \( e \)  is defined as $\xi_G(e) = d_G(u) + d_G(v) - 2$. The $minimum$ $edge$-$degree$ of \( G \), denoted by \( \xi(G) \), is defined as $ \min\left\{ \xi_G(e) \mid e \in E(G) \right\}$. A graph is simple if it has no parallel edges and loops.  Unless specified, all graphs in this paper refer to simple graphs.
For terminology and notation not defined here, we follow those provided in \cite{Bondy}.

If there exists a path between every pair of distinct vertices in $G$, then  $G$ is said to be  $connected$;  otherwise, $G$ is $disconnected$. Each maximal connected subgraph of $G$ is called its $component$. An edge subset \( S \subseteq E(G) \) is called an $edge$-$cut$ if \( G-S \) is disconnected. The $edge$-$connectivity$ \( \lambda(G) \) of $G$ is the minimum cardinality among all edge-cuts. It is well-known that \( \lambda(G) \leq \delta(G) \). So  the graph \( G \) is called $maximally$ $edge$-$connected$  if \( \lambda(G) = \delta(G) \); and the graph \( G \) is called $super$ $edge$-$connected$ if every minimum edge-cut isolates one vertex.

The edge-connectivity is one of the most fundamental parameters used to measure the reliability of networks. However, one deficiency of the edge-connectivity is that it only considers whether the remaining graph is disconnected, without taking into account the properties of the components. To address this limitation, Esfahanian and Hakimi \cite{Esfahanian} introduced the notion of
restricted edge-connectivity, which is a kind of conditional edge-connectivity initially proposed by Harary \cite{Harary}. If an edge subset \( S \subseteq E(G) \) satisfies that \( G - S \) is disconnected and each component of \( G - S \) has at least two vertices, then \( S \) is called a $restricted$ $edge$-$cut$. The $restricted$ $edge$-$connectivity$ of  \( G \), denoted by \( \lambda_2(G) \), is defined as the cardinality of a minimum restricted edge-cut in \( G \) if it has at least one; otherwise, \( \lambda_2(G) = +\infty \). In 1988, Esfahanian and Hakimi \cite{Esfahanian} proved that if \( G \) is not a star graph \( K_{1,n-1} \) and $G$ has at least four vertices, then \( \lambda_2(G) \leq \xi(G) \). So \( G \) is called $maximally$ $restricted$ $edge$-$connected$ if \( \lambda_2(G) = \xi(G) \); and \( G \) is called $super$ $restricted$ $edge$-$connected$ if every minimum restricted edge-cut isolates an edge.

Later on, the concept of \( k \)-restricted edge-connectivity was proposed by F{\`{a}}brega and Fiol \cite{Fabrega}. If an edge subset $S\subseteq E(G)$ satisfies that $G-S$ is disconnected and each component of $G-S$ has at least $k$ vertices, then $S$ is called a \( k \)-$restricted$ $edge$-$cut$. The $k$-$restricted$ $edge$-$connectivity$ $\lambda_k(G)$ is defined as the minimum cardinality among all $k$-restricted edge-cuts in $G$ if there exists one; otherwise, $\lambda_k(G) = +\infty$. From this definition, it follows that if \( \lambda_{k}(G) \) is finite, then \( \lambda_{l}(G) \) is finite for all integers $l$ satisfying $1\leq l\leq k$, and the inequalities $\lambda_{1}(G)
\leq\cdots\leq\lambda_{l}(G)\leq\cdots\leq\lambda_{k}(G)$ hold. Clearly, \( \lambda_1(G) \) is the edge-connectivity and \( \lambda_2(G) \) is the restricted edge-connectivity of $G$.

For two nonempty subsets \( X, Y \subseteq V(G) \), \( [X, Y]_G \) denotes the set of edges with one end in \( X \) and the other in \( Y \). When \( Y = V(G) \setminus X \), the edge set \( [X, Y]_G \) is simply written as \( \partial_G(X) \). Let \( \xi_3(G)= \min \{ |\partial_G(X)|\colon X \subseteq V(G), |X| = 3 \text{ and } G[X] \text{ is connected} \} \). Wang and Li \cite{Wang} provided necessary and sufficient conditions for the existence of 3-restricted edge-cut in graphs and showed that \( \lambda_3(G) \leq \xi_3(G) \) holds when there is at least one 3-restricted edge-cut in $G$. So  \( G \) is said to be $maximally$ $3$-$restricted$ $edge$-$connected$ if \( \lambda_3(G) = \xi_3(G) \), and  \( G \) is said to be $super$ $3$-$restricted$ $edge$-$connected$ if every minimum 3-restricted edge-cut isolates a component of order three.

Let \( G \) and \( H \) be two graphs. The \textit{direct product} (also known as the \textit{Kronecker product}, \textit{tensor product}, or \textit{cross product}) \( G \times H \) has the vertex set \( V(G \times H) = V(G) \times V(H) \), where two vertices \( (u_1, v_1) \) and \( (u_2, v_2) \) are adjacent in \( G \times H \) if and only if \( u_1u_2 \in E(G) \) and \( v_1v_2 \in E(H) \). Weichsel \cite{Weichsel} proved that the direct product of two nontrivial graphs is connected if and only if both graphs are connected and at least one of them is non-bipartite.

Bre{\v{s}}ar and {\v{S}}pacapan \cite{Bresar} gave some bounds on the edge-connectivity of the direct product of graphs.
Later, Cao, Brglez, {\v{S}}pacapan and Vumar \cite{Cao} obtained the edge-connectivity of the direct product of a nontrivial graph and a complete graph. {\v{S}}pacapan  \cite{Spacapan} not only obtained the edge-connectivity of the direct product of two graphs but also characterized the structure of the minimum edge-cut for the direct product of two graphs.

The $path$, the $cycle$, and the $complete$ $graph$ on $n$ vertices are denoted by $P_n$, $C_n$ and $K_n$, respectively. The $total$ $graph$ \( T_n \) is constructed from $K_n$ by attaching a loop to every vertex of \( K_n \). We use $K_{s,t}$ to denote the $complete$ $bipartite$ $graph$ with one part has $s$ vertices and the other part has $t$ vertices.   Ma, Wang and Zhang \cite{Ma} studied the restricted edge-connectivity of the direct product of a nontrivial graph with a complete graph or a total graph.
\begin{theorem}(\cite{Ma}\label{1})
For any nontrivial connected graph $G$ and any integer $n\geq 3$,$$\lambda_2\left(G\times K_n\right)=\min\left\{(n-1)\xi(G)+2(n-2),n(n-1)\lambda_2(G)\right\}.$$
\end{theorem}

\begin{theorem}(\cite{Ma}\label{2})
For any nontrivial connected graph $G$ and any integer $n\geq 3$,$$\lambda_2(G\times T_n)=min\{n\xi(G)+2(n-1),n^2\lambda_2(G)\}.$$
\end{theorem}

Bai, Tian and Yin \cite{Bai} further gave sufficient conditions for the direct product of a nontrivial graph with a complete graph or a total graph to be super restricted edge-connected.

\begin{theorem}(\cite{Bai}\label{3})
For any nontrivial connected graph $G$ and any integer $n \geq 5$. If $ n(n-1)\lambda_2 (G) > (n-1)\xi (G) + 2(n-2) $, then $ G \times K_{n} $ is super restricted edge-connected.
\end{theorem}

\begin{theorem}(\cite{Bai}\label{4})
For any nontrivial graph $G$ and any integer $n \geq 3$. If $n^{2}\lambda_2(G) > n\xi(G) + 2(n - 1)$, then $G \times T_n$ is super restricted edge-connected.
\end{theorem}

Guo, Hu, Yang and Zhao \cite{Guo} obtained the restricted edge-connectivity of the direct product of a connected graph $G$ satisfying $|V(G)| \leq n$ or $\Delta(G) \leq n-1$ with an odd cycle $C_n$.

\begin{theorem}(\cite{Guo}\label{5})
For any connected graph $G$ with $|V(G)|\leq n$ or $\Delta(G)\leq n-1$,$$\lambda_2 (G \times C_n) = \min \left\{ 2n\lambda_2(G), \min_{xy \in E(G)} 2(d_G(x) + d_G(y))- 2 \right\},$$ where $n\geq3$ and $n$ is odd.
\end{theorem}

Inspired by the aforementioned results, this paper investigates the 3-restricted edge-connectivity of the direct product of a $k$-regular connected graph with a cycle, a complete graph and a total graph. As corollaries, we establish sufficient conditions for these direct product graphs to be maximally 3-restricted edge-connected.\@ The subsequent section will introduce relevant definitions and lemmas. The main results will be presented in Section 3. The conclusion remarks are given in the last section.

\section{Preliminary}

Let \( \mathcal{G} = G \times H \). For any vertex \( u \in V(G) \), define \( H^u = \{(u, v) \in V(\mathcal{G}):v \in V(H)\} \), referred to as the \( H \)-layer of \( u \) in \( \mathcal{G} \). By the definition of the direct product, \( H^u \) is an independent set in \( \mathcal{G} \), and \( [H^{u_1}, H^{u_2}]_\mathcal{G} \neq \emptyset \) if and only if \( u_1u_2 \in E(G) \). For any vertex subset \( S \subseteq V(G) \), we define \( H^S = \bigcup_{u \in S} H^u \).

In the following two lemmas, Wang and Li \cite{Wang} established a necessary and sufficient condition for the existence of 3-restricted edge-cut and provided an upper bound for the 3-restricted edge-connectivity of $G$.

\begin{lemma}(\cite{Wang}\label{3})
The graph $G$ has a 3-restricted edge-cut if and only if $G$ has two vertex-disjoint paths of order 3.
\end{lemma}

\begin{lemma}(\cite{Wang}\label{4})
Let \( G \) be a simple connected graph of order at least 6. If $G$ has a 3-restricted edge-cut , then \( \lambda_{3}(G) \leq \xi_3(G) \).
\end{lemma}

\begin{lemma}\label{5}
Let \( \mathcal{G} = G \times H \), where \(G\) and \(H\) are two connected graphs. Let \(S\) be the minimum 3-restricted edge-cut of \( \mathcal{G}\), and let \(D_1\) and \(D_2\) be the two components of \( \mathcal{G} - S\). If there exist two non-adjacent edges \(u_1u_2\) and \(u_3u_4\) in \(G\) such that \(H^{u_1} \cup H^{u_2} \subseteq V(D_1)\) and \(H^{u_3} \cup H^{u_4} \subseteq V(D_2)\), then \(|S| \geq 2e(H)\lambda_2(G)\).
\end{lemma}

\noindent{\bf Proof.} Since there exist two non-adjacent edges \(u_1u_2\) and \(u_3u_4\) in \(G\), we obtain that \(G\) has a restricted edge-cut and $\lambda_2(G)$ is finite. There are at least \(\lambda_{2}(G)\) edge-disjoint paths from \(\{u_1, u_2\}\) to \(\{u_3, u_4\}\) in \(G\). Thus, there are at least \(2e(H)\lambda_2(G)\) edge-disjoint paths from \(H^{u_1} \cup H^{u_2}\) to \(H^{u_3} \cup H^{u_4}\) in \(\mathcal{G}\). Therefore, we have \(|S| \geq 2e(H)\lambda_{2}(G)\). $\Box$

The following five lemmas are essential for establishing the main results in the subsequent section.

\begin{lemma}(\cite{Bai}\label{6})
Let \( A \) be a vertex subset of the graph \( \mathcal{G} = K_2 \times K_n \) $(n \geq 5)$. Assume \( 2 \leq |A| \leq 2n - 2\). Then \( |[A, V(\mathcal{G}) \setminus A]_{\mathcal{G}}| \geq 2(n - 2) \), and the equality holds if and only if (i) \( |A| = 2 \) and \( \mathcal{G}[A] \) is isomorphic to \( K_2 \), or (ii) \( |A| = 2n - 2 \) and \( \mathcal{G}[V(\mathcal{G}) \setminus A] \) is isomorphic to \( K_2 \).
\end{lemma}

\begin{lemma}(\cite{Bai}\label{7})
Let  $A$  be a vertex subset of the graph  $\mathcal{G} = K_{2} \times T_{n}$  $(n \geq 3)$. Assume  $2 \leq |A| \leq 2n - 2$. Then \( |[A, V(\mathcal{G}) \setminus A]_{\mathcal{G}}| \geq 2(n-1) \), and the equality holds if and only if (i)  $|A| = 2$  and  $\mathcal{G}[A]$  is isomorphic to  $K_{2}$ , or (ii) $|A|
 = 2n - 2$  and \( \mathcal{G}[V(\mathcal{G}) \setminus A] \) is isomorphic to  $K_{2}$ .
\end{lemma}

\begin{lemma}\label{8}
Let \( B \) be a vertex subset of the graph \( \mathcal{G} = K_2 \times K_n \) $( n \geq 5)$. Assume \( 3 \leq |B| \leq 2n - 3 \). Then \( |[B, V(\mathcal{G}) \setminus B]_{\mathcal{G}}| \geq 3n - 7 \), and the equality holds if and only if (i) \( |B| = 3 \) and \( \mathcal{G}[B] \) is isomorphic to \( P_3 \), or (ii) \( |B| = 2n - 3 \) and \( \mathcal{G}[V(\mathcal{G}) \setminus B] \) is isomorphic to \( P_3 \), or (iii) $n=5$, \( \mathcal{G}[B] \) and \( \mathcal{G}[V(\mathcal{G}) \setminus B] \) are both isomorphic to \( K_{2,3} \).
\end{lemma}

\noindent{\bf Proof.} Let $U_{1}=\{u_{1}\}\times V(K_{n})$ and $U_{2}=\{u_{2}\}\times V(K_{n})$ be the bipartition of $\mathcal{G}$, where $\{u_1,u_2\}=V(K_2).$ Assume $B_i=B\cap U_i$ and $b_i=|B_i|$ for $i=1,2$. Without loss of generality, assume $|B|\leq|V(\mathcal{G})\setminus B|$, that is $|B|\leq n$. Then \( 3 \leq |B|=b_1 + b_2 \leq n \).

If \(b_1 = 0\) or \(b_2 = 0\), then $\left| [B, V(\mathcal{G}) \setminus B]_{\mathcal{G}} \right|=|B|(n - 1) \geq 3n - 3.$

If \(b_1>0\) and \(b_2>0\), then
\[
\begin{aligned}
\left| [B, V(\mathcal{G}) \setminus B]_{\mathcal{G}} \right| &\geq b_1(n - 1 - b_2) + b_2(n - 1 - b_1) \\
&= (n - 1)(b_1 + b_2) - 2b_1b_2 \\
&= (n - 1)|B| - 2b_1b_2.
\end{aligned}
\]

Since \( b_1 + b_2 = |B| \) and \( b_1, b_2 \) are positive integers, we have \( (n - 1)|B| - 2b_1b_2 \geq (n - 1)|B| - 2\left\lceil \frac{|B|}{2} \right\rceil \left\lfloor \frac{|B|}{2} \right\rfloor \). When \( |B| \) is an odd integer, \( (n - 1)|B| - 2\left\lceil \frac{|B|}{2} \right\rceil \left\lfloor \frac{|B|}{2} \right\rfloor = (n - 1)|B| - 2\frac{|B|+1}{2}\frac{|B|-1}{2} =-\frac{|B|^2}{2} + (n - 1)|B| + \frac{1}{2} \). When \( |B| \) is an even integer, \( (n - 1)|B| - 2\left\lceil \frac{|B|}{2} \right\rceil \left\lfloor \frac{|B|}{2} \right\rfloor = - \frac{|B|^2}{2}+(n - 1)|B|  \).

If \( |B| \) is an odd integer, then $\left| [B, V(\mathcal{G}) \setminus B]_{\mathcal{G}} \right|\geq -\frac{|B|^2}{2} + (n - 1)|B| + \frac{1}{2}\geq\min\{-\frac{3^2}{2} + 3(n - 1)+\frac{1}{2}, -\frac{n^2}{2}+n(n-1)+ \frac{1}{2}\}=\min\{3n - 7, \frac{n^2}{2} - n + \frac{1}{2}\}$ by $3\leq|B|\leq n$.
For \( n \geq 6 \), we have \( \frac{n^2}{2} - n + \frac{1}{2} > 3n - 7 \). For \( n = 5 \), we have \( \frac{n^2}{2} - n + \frac{1}{2} = 8 = 3n - 7 \). Thus we conclude that \( |[B, V(\mathcal{G}) \setminus B]_{\mathcal{G}}| \geq 3n - 7 \), equality holds if  $|B|=3$ and \(\mathcal{G}[B] \cong P_3\), or $n=5, |B|=5$, \( \mathcal{G}[B] \) and \( \mathcal{G}[V(\mathcal{G}) \setminus B] \) are both isomorphic to \( K_{2,3} \).

If \( |B| \) is an even integer, then $\left| [B, V(\mathcal{G}) \setminus B]_{\mathcal{G}} \right|\geq -\frac{|B|^2}{2} + (n - 1)|B| \geq\min\{-\frac{4^2}{2} + 4(n - 1), -\frac{n^2}{2}+n(n-1)\}=\min\{4n-12, \frac{n^2}{2} - n\}$ by $4\leq|B|\leq n$.
Since $\min\{4n-12, \frac{n^2}{2} - n\}>3n-7$, we have \( |[B, V(\mathcal{G}) \setminus B]_{\mathcal{G}}|>3n - 7 \) in this case.
$\Box$

By a similar argument as Lemma 2.5, we have the following lemma.

\begin{lemma}\label{9}
Let \( B \) be a vertex subset of the graph \( \mathcal{G} = K_2 \times T_n \) $( n \geq 3)$. Assume \( 3 \leq |B| \leq 2n - 3 \). Then \( |[B, V(\mathcal{G}) \setminus B]_{\mathcal{G}}| \geq 3n - 4 \), and the equality holds if and only if (i) \( |B| = 3 \) and \( \mathcal{G}[B] \) is isomorphic to \( P_3 \), or (ii) \( |B| = 2n - 3 \) and \( \mathcal{G}[V(\mathcal{G}) \setminus B] \) is isomorphic to \( P_3 \).
\end{lemma}

\begin{lemma}\label{10}
($i$) If \( n \) is an odd integer and \( n \geq 3 \), then \(\lambda(K_2 \times C_n) = \lambda_2(K_2 \times C_n) = \lambda_3(K_2 \times C_n) = 2\).

($ii$) If \( n \geq 5 \), then \(\lambda(K_2 \times K_n) = n-1\), \(\lambda_2(K_2 \times K_n) = 2n-4\) and \(\lambda_3(K_2 \times K_n) = 3n-7\).

($iii$) If \( n \geq 5 \), then \(\lambda(K_2 \times T_n) = n\), \(\lambda_2(K_2 \times T_n) = 2n-2\) and \(\lambda_3(K_2 \times T_n) = 3n-4\).
\end{lemma}

\noindent{\bf Proof.} ($i$) Since \(K_2 \times C_n\) is isomorphic to \(C_{2n}\) for an odd integer \(n\), it follows that \(\lambda(K_2 \times C_n) = \lambda_2(K_2 \times C_n) = \lambda_3(K_2 \times C_n) = 2\).

($ii$) Let $S$ be a minimum edge-cut of $K_2 \times K_n$. If $K_2 \times K_n-S$ contains an isolated vertex, then $|S|\geq n-1$. Otherwise, by Lemma 2.4, we have $|S|\geq 2n-4\geq n-1$. Thus
 \(\lambda(K_2 \times K_n)=|S| \geq n-1\). On the other hand, \(\lambda(K_2 \times K_n) \leq\delta(K_2 \times K_n) = n-1\). Therefore, it follows that \(\lambda(K_2 \times K_n) = n-1\).

Since \(K_2 \times K_n\) is not isomorphic to the star graph \(K_{1,n-1}\) and contains at least four vertices, we have \(\lambda_2(K_2 \times K_n) \leq \xi_2(K_2 \times K_n) = 2n-4\).  Lemma 2.4 implies \(\lambda_2(K_2 \times K_n) \geq 2n-4\). Therefore, \(\lambda_2(K_2 \times K_n) = 2n-4\).

By Lemmas 2.1 and 2.2,  we know that   \(K_2 \times K_n\) has a 3-restricted edge-cut and \(\lambda_3(K_2 \times K_n) \leq \xi_3(K_2 \times K_n) = 3n - 7\). Furthermore, by Lemma 2.5, it follows that \(\lambda_3(K_2 \times K_n) \geq 3n - 7\). Thus, \(\lambda_3(K_2 \times K_n) = 3n - 7\).

($iii$) The proof is similar to ($ii$).
$\Box$

\section{Main Results}

The $girth$ $g(G)$ of $G$ is the length of its shortest cycle if $G$ contains cycles; otherwise, define $g(G) = +\infty$. The only 1-regular  connected graph is $K_2$. By Lemma 2.7 ($i$), we have $\lambda_3(K_2 \times C_n) = 2$ for an odd integer  $n$. So we consider \(k \geq 2\) in the following theorem.

\begin{theorem}\label{4}
Let \( G \) be a $k$ $(\geq 2)$-regular connected graph with at least four vertices. Then
    \[
    \lambda_{3}(G \times C_{n}) = \begin{cases}
        \min\{2n\lambda_{2}(G), 6k - 6\},& \text{if } g(G) = 3; \\
        \min\{2n\lambda_{2}(G), 6k - 4\},& \text{if } g(G) \geq 4,
    \end{cases}
    \]
    where \( n \geq 3 \) and \( n \) is odd.
\end{theorem}
\noindent{\bf Proof.} Let \(\mathcal{G} = G \times C_{n}\). By Lemma 2.1, \(\mathcal{G}\) has a 3-restricted edge-cut. When \(g(G) = 3\), the definition of the direct product graph implies that \(g(\mathcal{G}) = 3\). According to Lemma 2.2, we have \(\lambda_{3}(\mathcal{G}) \leq \xi_{3}(\mathcal{G}) = 6k - 6\). When \(g(G) \geq 4\), the definition of the direct product graph implies that \(g(\mathcal{G}) \geq 4\). By Lemma 2.2, we obtain \(\lambda_{3}(\mathcal{G}) \leq \xi_{3}(\mathcal{G}) = 6k - 4\). Since $G$ is not a star graph and has at least four vertices, it follows $G$ has a restricted edge-cut. Let \(R\) be the minimum restricted edge-cut of the graph \(G\). Then \(G - R\) has exactly two nontrivial components \(X_1\) and \(X_2\). Define \(Y_i = X_i \times V(C_n)\) for \(i = 1, 2\). Since both \(\mathcal{G}[Y_1]\) and \(\mathcal{G}[Y_2]\) contain at least three vertices, it follows that \(\lambda_3(\mathcal{G}) \leq |[Y_1, Y_2]_\mathcal{G}| = 2n\lambda_{2}(G)\). Therefore,
\[
\lambda_3(G \times C_n) \leq
\begin{cases}
\min\{2n\lambda_2(G),\, 6k - 6\}, & \text{if } g(G) = 3; \\
\min\{2n\lambda_2(G),\, 6k - 4\}, & \text{if } g(G) \geq 4.
\end{cases}
\]

Now, it suffices to prove that
$$\lambda_3(G \times C_n) \geq \begin{cases}
\min\{2n\lambda_2(G), 6k - 6\}, & \text{if } g(G) = 3; \\
\min\{2n\lambda_2(G), 6k - 4\}, & \text{if } g(G) \geq 4.
\end{cases}$$

Let \( S \) be the minimum 3-restricted edge-cut of \( \mathcal{G} \). Then \( \mathcal{G} - S \) has exactly two components, say \( D_1 \) and \( D_2 \), where \( |V(D_1)| \geq 3 \) and \( |V(D_2)| \geq 3 \).

If for any vertex $u \in V(G)$, we have $C_n^{u} \subseteq V(D_1)$ or $C_n^{u} \subseteq V(D_2)$.\@ Considering that \( D_1 \) is connected and \( C_n^{u} \) is an independent set, it follows that there exist at least two layers \( C_n^{u_1} \) and \( C_n^{u_2} \) adjacent in \( D_1 \).\@ Similarly, there are at least two layers \( C_n^{u_3} \) and \( C_n^{u_4} \) adjacent
in \( D_2 \). By Lemma 2.3, we have $|S| \geq 2n\lambda_{2}(G).$
So we assume that  there exists a vertex \( u_1 \in V(G) \) such that \( C_n^{u_1} \cap V(D_1) \neq \emptyset \) and \( C_n^{u_1} \cap V(D_2) \neq \emptyset \) in the following proof. By Lemma 2.7, we have $\lambda(K_2 \times C_n) = 2$. Thus, for any vertex $u \in N_G(u_1)$, it follows that $\left| [C_n^{u_1}, C_n^u]_{\mathcal{G}} \cap S \right| \geq 2$.

\noindent{\bf Case 1.} Every vertex \( u \in N_{G}(u_1) \) satisfies \( C_n^{u} \subseteq V(D_1) \) or \( C_n^{u} \subseteq V(D_2) \).

Let \( M_1 = \{x_1, \ldots, x_{k_1}\} \) be the neighbors of \( u_1 \) in \( G \) such that \( C_n^{x_i} \subseteq V(D_1) \) for \( i \in \{1, \ldots, k_1\} \), and let \( M_2 = \{y_1, \ldots, y_{k_2}\} \) be the neighbors of \( u_1 \) in \( G \) such that \( C_n^{y_j} \subseteq V(D_2) \) for \( j \in \{1, \ldots, k_2\} \), where \( k_1 + k_2 = k \). Since \(C_n^{u_1}\) is an independent set, there exists at least one vertex \(x_1 \in M_1\) such that \(C_n^{x_1} \subseteq V(D_1)\) and one vertex \(y_1 \in  M_2 \) such that \(C_n^{y_1} \subseteq V(D_2)\). Therefore, \( |[C_n^{u_1}, C_n^{x_1}]_{\mathcal{G}} \cap S| + |[C_n^{u_1}, C_n^{y_1}]_{\mathcal{G}} \cap S| = 2n \).

\noindent{\bf Subcase 1.1.}\@ There exists a vertex \( x' \in N_{G}(M_1) \) and a vertex \( y' \in N_{G}(M_2) \) such that \( C_n^{x'} \subseteq V(D_1) \) and \( C_n^{y'} \subseteq V(D_2) \).

Without loss of generality, assume \( x' \in N_{G}(x_1) \) and \( y' \in N_{G}(y_1) \). Then \( C_n^{x'} \cup C_n^{x_1} \subseteq V(D_1) \) and \( C_n^{y'} \cup C_n^{y_1} \subseteq V(D_2) \). By Lemma 2.3, we have \( |S| \geq 2n\lambda_2(G) \).

\noindent{\bf Subcase 1.2.} Each vertex \( x \in N_{G}(M_1) \) satisfies \( C_n^x \cap V(D_2) \neq \emptyset \), or each vertex \( y \in N_{G}(M_2) \) satisfies \( C_n^y \cap V(D_1) \neq \emptyset \).

Without loss of generality, assume that each vertex \( x \in N_G(M_1) \) satisfies \( C_n^x \cap V(D_2) \neq \emptyset \). Consequently, for any vertex \( x \in N_G(x_1) \), we have \( C_n^x \cap V(D_2) \neq \emptyset \). By  \(\lambda(K_2 \times C_n) = 2\), it follows that $\left| \left[ C_n^{x_1}, C_n^x \right]_{\mathcal{G}} \cap S \right| \geq 2$ for any vertex \( x \in N_G(x_1) \).

If for every vertex \( x \in N_{G}(x_1)\setminus\{u_1\} \), it holds that \( C_n^x \subseteq V(D_2) \), then \( \bigl|[C_n^{x_1}, C_n^{x}]_{\mathcal{G}} \cap S\bigr| = 2n \). Therefore,
$$\begin{aligned}
\left|S\right| & \geq \left|[C_{n}^{u_1},C_{n}^{x_1}]_{\mathcal{G}}\cap S\right| + \left|[C_{n}^{u_1},C_{n}^{y_1}]_{\mathcal{G}}\cap S\right| +\sum_{u\in N_{G}(u_1)\setminus\{x_1,y_1\}} \left|[C_{n}^{u_1},C_{n}^{u}]_{\mathcal{G}}\cap S\right| \\
& + \sum_{x\in N_{G}(x_1)\setminus\{u_1\}} \left|[C_{n}^{x_1},C_{n}^{x}]_{\mathcal{G}}\cap S\right| \\
& \geq 2n + 2(k-2) + 2n(k-1) \\
& = 2nk + 2k - 4 \\
& > 6k - 4.
\end{aligned}$$

If there exists a vertex \( x' \in N_{G}(x_1) \setminus \{u_1\} \) such that \( C_n^{x'} \cap V(D_1) \neq \emptyset \) and \( C_n^{x'} \cap V(D_2) \neq \emptyset \), then \( x' \notin N_{G}(u_1) \). By   \(\lambda(K_2 \times C_n) = 2\), it follows that $|[ C_n^{x'}, C_n^x]_{\mathcal{G}} \cap S | \geq 2$  for any vertex \( x \in N_G(x') \). Therefore,
\begin{align*}
|S| &\geq \left|[C_n^{u_1}, C_n^{x_1}]_{\mathcal{G}} \cap S\right| + \left|[C_n^{u_1}, C_n^{y_1}]_{\mathcal{G}} \cap S\right|
+ \sum_{\substack{u \in N_{G}(u_1)  \setminus \{x_1, y_1\}}} \left|[C_n^{u_1}, C_n^{u}]_{\mathcal{G}} \cap S\right| \\
&+ \sum_{\substack{x \in N_{G}(x_1) \setminus \{u_1\}}} \left|[C_n^{x_1}, C_n^{x}]_{\mathcal{G}} \cap S\right|
+ \sum_{\substack{x \in N_{G}(x') \setminus \{x_1\}}} |[C_n^{x'}, C_n^{x}]_{\mathcal{G}} \cap S| \\
&= 2n + 2(k-2) + 2(k-1) + 2(k-1) \\
&= 6k - 4 + 2n - 4 \\
&> \mathrm{6k - 4}.
\end{align*}

\noindent{\bf Case 2.} There is a vertex \( u_2 \in N_{G}(u_1) \) such that \( C_n^{u_2} \cap V(D_1) \neq \emptyset \) and \( C_n^{u_2} \cap V(D_2) \neq \emptyset \).

Let $A_1 = C_n^{u_1} \cap V(D_1), B_1 = C_n^{u_1} \cap V(D_2), A_2 = C_n^{u_2} \cap V(D_1), B_2 = C_n^{u_2} \cap V(D_2).$ Then \( |A_1| + |B_1| = |A_2| + |B_2| = n \) and \( |A_i| \geq 1 \), \( |B_i| \geq 1 \) for \( i = 1, 2 \).
By  $\lambda(K_2 \times C_n) =2$, we have $\left| \left[ C_n^{u_2}, C_n^u \right]_{\mathcal{G}} \cap S \right| \geq 2$ for any vertex $u \in N_G(u_2)$.

\noindent{\bf Subcase 2.1.}\@ There is a vertex \( u_3 \in N_{G}(u_1) \setminus \{u_2\} \) such that \( C_n^{u_3} \cap V(D_1) \neq \emptyset \) and \( C_n^{u_3} \cap V(D_2) \neq \emptyset \).

By   $\lambda(K_2 \times C_n) = 2$, it follows $\left| \left[ C_n^{u_3}, C_n^u \right]_{\mathcal{G}} \cap S \right| \geq 2$ for any vertex $u \in N_G(u_3)$.

If \( u_3 \notin N_{G}(u_2) \), then
\begin{align*}
	|S| &\geq \sum_{\substack{u \in N_{G}(u_1) }} |[C_n^{u_1},C_n^u]_\mathcal{G} \cap S| +\sum_{\substack{u \in N_{G}(u_2) \backslash \{u_1\}}} |[C_n^{u_2},C_n^u]_\mathcal{G} \cap S| \\
	& +\sum_{\substack{u \in N_{G}(u_3) \backslash \{u_1\}}} |[C_n^{u_3},C_n^u]_\mathcal{G} \cap S| \\
	&\geq 2k+ 2(k-1) + 2(k-1) \\
	&= 6k-4.
\end{align*}

If \( u_3 \in N_{G}(u_2) \), then
\begin{align*}
	|S| &\geq \sum_{\substack{u \in N_{G}(u_1) }} |[C_n^{u_1},C_n^u]_\mathcal{G} \cap S| +\sum_{\substack{u \in N_{G}(u_2) \backslash \{u_1\}}} |[C_n^{u_2},C_n^u]_\mathcal{G} \cap S| \\
	& +\sum_{\substack{u \in N_{G}(u_3) \backslash \{u_1,u_2\}}} |[C_n^{u_3},C_n^u]_\mathcal{G} \cap S| \\
	&\geq 2k+ 2(k-1) + 2(k-2) \\
	&= 6k-6.
\end{align*}

\noindent{\bf Subcase 2.2.} For every vertex \( u \in N_{G}(u_1) \setminus \{u_2\} \), either \( C_n^{u} \subseteq V(D_1) \) or \( C_n^{u} \subseteq V(D_2) \).

If there exists a vertex \( u_3' \in N_G(u_2) \setminus \{u_1\} \) such that \( C_n^{u_3'} \cap V(D_1) \neq \emptyset \) and \( C_n^{u_3'} \cap V(D_2) \neq \emptyset \), then, analogously to Subcase 2.1, we obtain either \(|S| \geq 6k-4\) or \(|S| \geq 6k-6\). So we assume that for any vertex \( u \in (N_G(u_1) \cup N_G(u_2)) \setminus \{u_1, u_2\} \), either \( C_n^u \subseteq V(D_1) \) or \( C_n^u \subseteq V(D_2) \).

\noindent{\bf Subcase 2.2.1.}\@ Any vertex \( u \in (N_G(u_1) \cup N_G(u_2)) \setminus \{u_1, u_2\} \) satisfies \( C_n^u \subseteq V(D_1) \), or any vertex \( u \in (N_G(u_1) \cup N_G(u_2)) \setminus \{u_1, u_2\} \) satisfies \( C_n^u \subseteq V(D_2) \).

Without loss of generality, assume that each vertex \( u \in (N_G(u_1) \cup N_G(u_2)) \setminus \{u_1, u_2\} \) satisfies \( C_n^u \subseteq V(D_1) \).  Then $V(D_{2}) \subseteq C_n^{u_1} \cup C_n^{u_2}$ and   $|B_1| + |B_2| \geq 3$. Thus $\sum_{u \in N_{G}(u_{1}) \backslash \{u_{2}\}} \bigl|[C_{n}^{u_{1}}, C_{n}^{u}]_{\mathcal{G}} \cap S\bigr| + \sum_{u \in N_{G}(u_{2}) \backslash \{u_{1}\}} \bigl|[C_{n}^{u_{2}}, C_{n}^{u}]_{\mathcal{G}} \cap S\bigr| =2(k - 1)|B_1|+2(k - 1)|B_2|\geq3 \times 2(k - 1) = 6(k - 1)$.

Therefore,
\[
\begin{aligned}
|S| &\geq \bigl|[C_n^{u_1}, C_n^{u_2}]_{\mathcal{G}} \cap S\bigr| + \sum_{u \in N_{G}(u_1) \backslash \{u_2\}} \bigl|[C_n^{u_1}, C_n^{u}]_{\mathcal{G}} \cap S\bigr| \\
&\quad + \sum_{u \in N_{G}(u_2) \backslash \{u_1\}} \bigl|[C_n^{u_2}, C_n^{u}]_{\mathcal{G}} \cap S\bigr| \\
&\geq 2 + 6(k - 1) \\
&= 6k - 4.
\end{aligned}
\]

\noindent{\bf Subcase 2.2.2.} There exists a vertex \(u_4 \in (N_G(u_1) \cup N_G(u_2)) \setminus \{u_1, u_2\}\) such that \(C_n^{u_4} \subseteq V(D_1)\), and there exists a vertex \(u_5 \in (N_G(u_1) \cup N_G(u_2)) \setminus \{u_1, u_2\}\) such that \(C_n^{u_5} \subseteq V(D_2)\).

If there exists a vertex \(u^{\prime} \in N_{G}(u_{4})\) such that \(C_{n}^{u^{\prime}} \subseteq V(D_{1})\), and there exists a vertex \(u^{\prime\prime} \in N_{G}(u_{5})\) such that \(C_{n}^{u^{\prime\prime}} \subseteq V(D_{2})\), then by Lemma 2.3, we have $|S| \geq 2n \lambda_{2}(G)$.
So we assume that
every vertex \(u \in N_{G}(u_{4})\) satisfies \(C_{n}^{u} \cap V(D_{2}) \neq \emptyset \), or every vertex \(u \in N_{G}(u_{5})\) satisfies \(C_{n}^{u} \cap V(D_{1}) \neq \emptyset \).  Without loss of generality, assume that every vertex \(u \in N_{G}(u_{4})\) satisfies \(C_{n}^{u} \cap V(D_{2}) \neq \emptyset \). By  $\lambda(K_2 \times C_n) = 2$ we have \(\left| [C_n^{u_4}, C_n^{u}]_{\mathcal{G}} \cap S \right| \geq 2\) for any vertex \( u \in N_G(u_4) \).

If \(u_4 \in N_{G}(u_1)\) and \(u_4 \in N_{G}(u_2)\), then

\[
\begin{aligned}
|S| &\geq \sum_{u \in N_{G}(u_1)} \bigl|[C_n^{u_1}, C_n^{u}]_\mathcal{G} \cap S\bigr| + \sum_{u \in N_G(u_2) \setminus \{u_1\}} \bigl|[C_n^{u_2}, C_n^{u}]_\mathcal{G} \cap S\bigr| \\
&\quad + \sum_{u \in N_G(u_4) \setminus \{u_1, u_2\}} \bigl|[C_n^{u_4}, C_n^{u}]_\mathcal{G} \cap S\bigr| \\
&\geq 2k + 2(k - 1) + 2(k - 2) \\
&= 6k - 6.
\end{aligned}
\]

If \(u_4 \in N_G(u_1)\) but \(u_4 \notin N_G(u_2)\), or \(u_4 \in N_G(u_2)\) but \(u_4 \notin N_G(u_1)\), then, without loss of generality, we assume \(u_4 \in N_G(u_1)\) and \(u_4 \notin N_G(u_2)\). Thus,

\[
\begin{aligned}
|S| &\geq \sum_{u \in N_{G}(u_1)} \bigl|[C_n^{u_1}, C_n^{u}]_\mathcal{G} \cap S\bigr| + \sum_{u \in N_G(u_2) \setminus \{u_1\}} \bigl|[C_n^{u_2}, C_n^{u}]_\mathcal{G} \cap S\bigr| \\
&\quad + \sum_{u \in N_G(u_4) \setminus \{u_1\}} \bigl|[C_n^{u_4}, C_n^{u}]_\mathcal{G} \cap S\bigr| \\
&\geq 2k + 2(k - 1) + 2(k - 1) \\
&= 6k - 4.
\end{aligned}
\]

The proof is thus complete. $\Box$

By Lemma 2.7, we have $\lambda_3(K_2 \times K_n) = 3n-7$. So we consider \(k \geq 2\) in the following theorem.

\begin{theorem}\label{3}
Let $G$ be a $k$ $(\geq 2)$-regular connected graph with at least four vertices. Then
$$\lambda_{3}(G\times K_{n})=\begin{cases}\min\{n(n - 1)\lambda_{2}(G),3k(n - 1)-6\},& \text{if } g(G) = 3;\\\min\{n(n - 1)\lambda_{2}(G),3k(n - 1)-4\},& \text{if } g(G)\geq4,\end{cases}$$
where $n\geq5$.
\end{theorem}
\noindent{\bf Proof.} Let \(\mathcal{G}=G \times K_{n}\). By Lemma 2.1, \(\mathcal{G}\) has a 3-restricted edge-cut. When \(g(G) = 3\), the definition of the direct product graph implies that \(g(\mathcal{G}) = 3\). According to Lemma 2.2, we have \(\lambda_{3}(\mathcal{G}) \leq \xi_{3}(\mathcal{G}) = 3k(n - 1) - 6\). When \(g(G) \geq 4\), the definition of the direct product graph implies that \(g(\mathcal{G}) \geq 4\). By Lemma 2.2, we obtain \(\lambda_{3}(\mathcal{G}) \leq \xi_{3}(\mathcal{G}) = 3k(n - 1) - 4\). Since $G$ is not a star graph and has at least four vertices, it follows $G$ has a restricted edge-cut. Let \(R\) be the minimum restricted edge-cut of the graph \(G\). Then \(G - R\) has exactly two nontrivial components \(X_1\) and \(X_2\). Define \(Y_i = X_i \times V(K_n)\) for \(i = 1, 2\). Since both \(\mathcal{G}[Y_1]\) and \(\mathcal{G}[Y_2]\) contain at least three vertices, it follows that \(\lambda_3(G) \leq |[Y_1, Y_2]_\mathcal{G}| = n(n - 1)\lambda_2(G)\). Therefore,
\[
\lambda_3(G \times K_n) \leq
\begin{cases}
\min\{n(n-1) \lambda_2(G), 3k(n-1) - 6\}, & \text{if }g(G) = 3;\\
\min\{n(n-1) \lambda_2(G), 3k(n-1) - 4\}, & \text{if } g(G) \geq 4.
\end{cases}
\]

Now, it suffices to prove that
\[
\lambda_3(G \times K_n) \geq
\begin{cases}
\min\{n(n-1) \lambda_2(G), 3k(n-1) - 6\}, & \text{if }g(G) = 3;\\
\min\{n(n-1) \lambda_2(G), 3k(n-1) - 4\}, & \text{if } g(G) \geq 4.
\end{cases}
\]

Let \( S \) be a minimum 3-restricted edge-cut of \( G \).\@ Then \( G - S \) has exactly two components, say \( D_1 \) and \( D_2 \), where \( |V(D_1)| \geq 3 \) and \( |V(D_2)| \geq 3 \).

If for any vertex \( u \in V(G) \), we have \( K_n^{u} \subseteq V(D_1) \) or \( K_n^{u} \subseteq V(D_2) \).\@ Considering that \( D_1 \) is connected and $K_n^{u}$ is an independent set, it follows that there exist at least two layers \( K_n^{u_1} \) and \( K_n^{u_2} \) in \( D_1 \) that are adjacent. Similarly, there are at least two layers \( K_n^{u_3} \) and \( K_n^{u_4} \) adjacent in $D_2$.\@  By Lemma 2.3, we have \( |S| \geq n(n-1) \lambda_2(G) \).
So we assume that  there exists a vertex \( u_1 \in V(G) \) such that \( K_n^{u_1} \cap V(D_1) \neq \emptyset \) and \( K_n^{u_1} \cap V(D_2) \neq \emptyset \) in the following proof. By Lemma 2.7, we have $\lambda(K_2 \times K_n) = n-1$. Thus, for any vertex $u \in N_G(u_1)$, it follows that $\left| [K_n^{u_1}, K_n^u]_{\mathcal{G}} \cap S \right| \geq n-1$.

\noindent{\bf Case 1.} Every vertex \( u \in N_{G}(u_1) \) satisfies \( K_n^{u} \subseteq V(D_1) \) or \( K_n^{u} \subseteq V(D_2) \).

Let \( M_1 = \{x_1, \dots, x_{k_1}\} \) be the neighbors of \( u_1 \) in \( G \) such that \( K_n^{x_i} \subseteq V(D_1) \) for \( i \in \{1, \dots, k_1\} \), and let \( M_2 = \{y_1, \dots, y_{k_2}\} \) be the neighbors of \( u_1 \) in \( G \) such that \( K_n^{y_j} \subseteq V(D_2) \) for \( j \in \{1, \dots, k_2\} \), where \( k_1 + k_2 = k \). Since \( K_n^{u_1} \) is an independent set, there exists at least one vertex \( x_1 \in  M_1  \) such that \( K_n^{x_1} \subseteq V(D_1) \) and one vertex \( y_1 \in  M_2 \) such that \( K_n^{y_1} \subseteq V(D_2) \). Therefore, $ |[K_n^{u_1}, K_n^{x_1}]_G \cap S| + |[K_n^{u_1}, K_n^{y_1}]_G \cap S| =n(n-1).$

\noindent{\bf Subcase 1.1.}\@ There is a vertex \( x^{\prime} \in N_{G}(M_1) \) such that \( K_n^{x^{\prime}} \subseteq V(D_1) \), and there is a vertex \( y^{\prime} \in N_{G}(M_2) \) such that \( K_n^{y^{\prime}} \subseteq V(D_2) \).

Without loss of generality, assume \( {x^{\prime}} \in N_{G}(x_1) \) and \( {y^{\prime}} \in N_{G}(y_1) \). Then \( K_n^x \cup K_n^{x_1} \subseteq V(D_1) \) and \( K_n^y \cup K_n^{y_1} \subseteq V(D_2) \).  By Lemma 2.3, we have \( |S| \geq n(n-1)\lambda_2(G) \).

\noindent{\bf Subcase 1.2.} Each vertex \( x \in N_{G}(M_1) \) satisfies \( K_n^x \cap V(D_2) \neq \emptyset \), or each vertex \( y \in N_{G}(M_2) \) satisfies \( K_n^x \cap V(D_1) \neq \emptyset \).

Without loss of generality, assume that \( K_n^x \cap V(D_2) \neq \emptyset \) for each vertex \( x \in N_G(M_1) \). Consequently, for any vertex \( x \in N_G(x_1) \), we have \( K_n^x \cap V(D_2) \neq \emptyset \). By  $\lambda(K_2 \times K_n) = n-1$, it follows that $\left| [K_n^{x_1}, K_n^u]_{\mathcal{G}} \cap S \right| \geq n-1$ for any vertex $x \in N_G(x_1)$.

If for every vertex \( x \in N_{G}(x_1) \setminus \{u_1\} \), it holds that \( K_n^x \subseteq V(D_2) \), then
$|[K_n^{x_1}, K_n^{x}]_\mathcal{G} \cap S|= n(n-1). $
Thus,
\begin{align*}
|S| &\geq \left|[K_n^{u_1}, K_n^{x_1}]_{\mathcal{G}} \cap S\right| + \left|[K_n^{u_1}, K_n^{y_1}]_{\mathcal{G}} \cap S\right|+\sum_{\substack{u \in N_{G}(u_1) \setminus \{x_1, y_1\}}} \left|[K_n^{u_1}, K_n^{u}]_{\mathcal{G}} \cap S\right| \\
&+ \sum_{\substack{x \in N_{G}(x_1) \setminus \{u_1\}}} \left|[K_n^{x_1}, K_n^{x}]_{\mathcal{G}} \cap S\right| \\
&= n(n-1) + (k-2)(n-1) + n(n-1)(k-1) \\
&= (n-1)(nk + k - 2) \\
&> 3k(n-1) - 4.
\end{align*}

If there exists a vertex \(x^{\prime} \in N_{G}(x_1) \backslash \{u_1\}\) such that \(K_n^{x^{\prime}} \cap V(D_1) \neq \emptyset\) and \(K_n^{x^{\prime}} \cap V(D_2) \neq \emptyset\), then \(x^{\prime} \notin N_{G}(u_{1})\). By  $\lambda(K_2 \times K_n) = n-1$, it follows that $| [K_n^{x^{\prime}}, K_n^x]_{\mathcal{G}} \cap S | \geq n-1$ for any vertex $x \in N_G(x^{\prime})$. Therefore,
\begin{align*}
|S| &\geq \left|[K_n^{u_1}, K_n^{x_1}]_{\mathcal{G}} \cap S\right|+\left|[K_n^{u_1}, K_n^{y_1}]_{\mathcal{G}} \cap S\right|+\sum_{\substack{u \in N_{G}(u_1) \setminus \{x_1, y_1\}}} \left|[K_n^{u_1}, K_n^{u}]_{\mathcal{G}} \cap S\right| \\
&+ \sum_{\substack{x \in N_{G}(x_1) \setminus \{u_1\}}} \left|[K_n^{x_1}, K_n^{x}]_{\mathcal{G}} \cap S\right|+\sum_{\substack{x \in N_{G}(x') \setminus \{x_1\}}} |[K_n^{x'}, K_n^{x}]_{\mathcal{G}} \cap S| \\
&= n(n-1)+(k-2)(n-1)+(k-1)(n-1)+(k-1)(n-1) \\
&= (n-1)(3k + n - 4) \\
&> 3k(n-1) - 4.
\end{align*}

\noindent{\bf Case 2.} There is a vertex \( u_2 \in N_{G}(u_1) \) such that \( K_n^{u_2} \cap V(D_1) \neq \emptyset \) and \( K_n^{u_2} \cap V(D_2) \neq \emptyset \).

Let \( A_1 = K_n^{u_1} \cap V(D_1) \), \( B_1 = K_n^{u_1} \cap V(D_2) \), \( A_2 = K_n^{u_2} \cap V(D_1) \) and \( B_2 = K_n^{u_2} \cap V(D_2) \). Then, \( |A_1| + |B_1| = |A_2| + |B_2| = n \) and \( |A_i| \geq 1 \) and \( |B_i| \geq 1 \) for \( i = 1, 2 \). Thus, we obtain $2\leq|B_1| + |B_2|\leq2n-2$.
By $\lambda(K_2 \times K_n) = n-1$, we have $| [K_n^{u_2}, K_n^u]_{\mathcal{G}} \cap S| \geq n-1$ for any vertex $u \in N_G(u_2) \setminus \{u_1\}$.
By $\lambda_2(K_2 \times K_n) =2(n-2)$, we have
$\left| [K_n^{u_1}, K_n^{u_2}]_{\mathcal{G}} \cap S \right| \geq 2(n-2)$.

\noindent{\bf Subcase 2.1.}\@ There is a vertex \( u_3 \in N_{G}(u_1) \setminus \{u_2\} \) such that \( K_n^{u_3} \cap V(D_1) \neq \emptyset \) and \( K_n^{u_3} \cap V(D_2) \neq \emptyset \).

By $\lambda(K_2 \times K_n) = n-1$, we have $| [K_n^{u_3}, K_n^u]_{\mathcal{G}} \cap S| \geq n-1$ for any vertex $u \in N_G(u_3) \setminus \{u_1\}$. We obtain $\left| [K_n^{u_1}, K_n^{u_3}]_{\mathcal{G}} \cap S \right| \geq 2(n-2)$ by $\lambda_2(K_2 \times K_n) =2(n-2)$.

If \( u_3 \notin N_{G}(u_2) \), then
\begin{align*}
|S| &\geq \left|[K_n^{u_1}, K_n^{u_2}]_{\mathcal{G}} \cap S\right| + \left|[K_n^{u_1}, K_n^{u_3}]_{\mathcal{G}} \cap S\right|+\sum_{\substack{u \in N_{G}(u_1) \setminus \{u_2, u_3\}}} \left|[K_n^{u_1}, K_n^{u}]_{\mathcal{G}} \cap S\right| \\
&+ \sum_{\substack{u \in N_{G}(u_2) \setminus \{u_1\}}} \left|[K_n^{u_2}, K_n^{u}]_{\mathcal{G}} \cap S\right|+\sum_{\substack{u \in N_{G}(u_3) \setminus \{u_1\}}} \left|[K_n^{u_3}, K_n^{u}]_{\mathcal{G}} \cap S\right| \\
&\geq 2(n-2) + 2(n-2) + (k-2)(n-1) + (k-1)(n-1) + (k-1)(n-1) \\
&= (n-1)(2+2+k-2+k-1+k-1) - 4\\
&= 3k(n-1) - 4.
\end{align*}

If \( u_3 \in N_{G}(u_2) \), then  $\left| [K_n^{u_2}, K_n^{u_3}]_{\mathcal{G}} \cap S \right| \geq 2(n-2)$
by $\lambda_2(K_2 \times K_n) =2(n-2)$. Thus

\begin{align*}
|S| &\geq \left|[K_n^{u_1}, K_n^{u_2}]_{\mathcal{G}} \cap S\right| + \left|[K_n^{u_1}, K_n^{u_3}]_{\mathcal{G}} \cap S\right| + \left|[K_n^{u_2}, K_n^{u_3}]_{\mathcal{G}} \cap S\right|\\
&+\sum_{{u \in N_{G}(u_1) \setminus \{u_2, u_3\}}} \left|[K_n^{u_1}, K_n^{u}]_{\mathcal{G}} \cap S\right|+\sum_{{u \in N_{G}(u_2) \setminus \{u_1, u_3\}}} \left|[K_n^{u_2}, K_n^{u}]_{\mathcal{G}} \cap S\right|\\
&+\sum_{{u \in N_{G}(u_3) \setminus \{u_1, u_2\}}} \left|[K_n^{u_3}, K_n^{u}]_{\mathcal{G}} \cap S\right| \\
&\geq 2(n-2) + 2(n-2) + 2(n-2) + (k-2)(n-1) + (k-2)(n-1) + (k-2)(n-1) \\
&= (n-1)(2+2+2+k-2+k-2+k-2) - 6\\
&= 3k(n-1) - 6.
\end{align*}

\noindent{\bf Subcase 2.2.} For each vertex \( u \in N_{G}(u_1) \setminus \{u_2\} \), either \( K_n^{u} \subseteq V(D_1) \) or \( K_n^{u} \subseteq V(D_2) \).

If there exists a vertex \( u_3' \in N_G(u_2) \setminus \{u_1\} \) such that \( C_n^{u_3'} \cap V(D_1) \neq \emptyset \) and \( C_n^{u_3'} \cap V(D_2) \neq \emptyset \), then, analogously to Subcase 2.1, we obtain either \(|S| \geq 3n(k-1)-4\) or \(|S| \geq 3n(k-1)-6\). So we assume that, for each \(u \in N_{G}(u_2) \backslash \{u_1\}\), either \(K_{n}^{u} \subseteq V(D_{1})\) or \(K_{n}^{u} \subseteq V(D_{2})\).

\noindent{\bf Subcase 2.2.1.}\@ Any vertex \(u \in (N_{G}(u_{1}) \cup N_{G}(u_{2})) \backslash \{u_1, u_2\}\) satisfies \(K_{n}^{u} \subseteq V(D_{1})\), or any vertex \(u \in (N_{G}(u_{1}) \cup N_{G}(u_{2})) \backslash \{u_1, u_2\}\) satisfies \(K_{n}^{u} \subseteq V(D_{2})\).

Without loss of generality, assume that for any vertex \(u \in (N_{G}(u_{1}) \cup N_{G}(u_{2})) \backslash \{u_1, u_2\}\) such that \( K_n^u \subseteq V(D_1) \). Then $V(D_{2}) \subseteq K_n^{u_1} \cup K_n^{u_2}$ and thus $|B_1| + |B_2| \geq 3$. By \( |A_1| + |A_2| \geq 2 \), it follows that $3 \leq |B_1| + |B_2| \leq 2n - 2.$

If \( 3 \leq |B_1| + |B_2| \leq 2n - 3 \),
then $\left|[K_n^{u_1}, K_n^{u_2}]_{\mathcal{G}} \cap S\right| \geq 3n-7$ by Lemma 2.5. And we have $\sum_{\substack{u \in N_{G}(u_1) \setminus \{u_2\}}} \left|[K_n^{u_1}, K_n^{u}]_{\mathcal{G}} \cap S\right| + \sum_{\substack{u \in N_{G}(u_2) \setminus \{u_1\}}} \left|[K_n^{u_2}, K_n^{u}]_{\mathcal{G}} \cap S\right|
= |B_1|(k-1)(n-1)+|B_2|(k-1)(n-1)
\geq 3(k-1)(n-1).$ Thus,
\begin{align*}
|S| &\geq \left|[K_n^{u_1}, K_n^{u_2}]_{\mathcal{G}} \cap S\right|+\sum_{\substack{u \in N_{G}(u_1) \setminus \{u_2\}}} \left|[K_n^{u_1}, K_n^{u}]_{\mathcal{G}} \cap S\right|+\sum_{\substack{u \in N_{G}(u_2) \setminus \{u_1\}}} \left|[K_n^{u_2}, K_n^{u}]_{\mathcal{G}} \cap S\right| \\
&\geq 3n-7 + 3(k-1)(n-1) \\
&= (n-1)(3+3k-3) - 4\\
&= 3k(n-1) - 4.
\end{align*}

When \( |B_1| + |B_2| = 2n-2 \), we have
$\sum_{\substack{u \in N_{G}(u_1) \setminus \{u_2\}}} \left|[K_n^{u_1}, K_n^{u}]_{\mathcal{G}} \cap S\right| + \sum_{\substack{u \in N_{G}(u_2) \setminus \{u_1\}}} \left|[K_n^{u_2}, K_n^{u}]_{\mathcal{G}} \cap S\right| =|B_1|(k-1)(n-1)+|B_2|(k-1)(n-1)=(2n-2)(k-1)(n-1).$
Since \(n\geq5\), we have
$$
\sum_{u\in N_{G}(u_{1})\setminus\{u_{2}\}} \left|[K_{n}^{u_{1}},K_{n}^{u}]_{\mathcal{G}}\cap S\right|+\sum_{u\in N_{G}(u_{2})\setminus\{u_{1}\}} \left|[K_{n}^{u_{2}},K_{n}^{u}]_{\mathcal{G}}\cap S\right|>4(k - 1)(n - 1).
$$
Thus,
\begin{align*}
|S| &\geq \left|[K_n^{u_1}, K_n^{u_2}]_{\mathcal{G}} \cap S\right|+\sum_{\substack{u \in N_{G}(u_1) \setminus \{u_2\}}} \left|[K_n^{u_1}, K_n^{u}]_{\mathcal{G}} \cap S\right|+\sum_{\substack{u \in N_{G}(u_2) \setminus \{u_1\}}} \left|[K_n^{u_2}, K_n^{u}]_{\mathcal{G}} \cap S\right| \\
&> 2(n-2)+ 4(k-1)(n-1) \\
&= (n-1)(4k - 2) - 2 \\
&> 3k(n-1) - 4.
\end{align*}

\noindent{\bf Subcase 2.2.2.}\@ For every vertex \(u \in N_{G}(u_1) \backslash \{u_2\}\), it holds that \(K_{n}^{u} \subseteq V(D_{1})\), and for every vertex \(u \in N_{G}(u_2) \backslash \{u_1\}\), it holds that \(K_{n}^{u} \subseteq V(D_{2})\).
Or, for every vertex \(u \in N_{G}(u_1) \backslash \{u_2\}\), it holds that \(K_{n}^{u} \subseteq V(D_{2})\), and for every vertex \(u \in N_{G}(u_2) \backslash \{u_1\}\), it holds that \(K_{n}^{u} \subseteq V(D_{1})\).

Without loss of generality, assume that for every vertex \(u \in N_{G}(u_1) \backslash \{u_2\}\), it holds that \(K_{n}^{u} \subseteq V(D_{1})\), and for every vertex \(u \in N_{G}(u_2) \backslash \{u_1\}\), it holds that \(K_{n}^{u} \subseteq V(D_{2})\).

Let \(M_3 = \{z_1, ..., z_{k-1}\}\) be the set of neighbors of \(u_1\) in \(G\) except for $u_{2}$ and let \(M_4 = \{w_1, ..., w_{k-1}\}\) be the set of neighbors of \(u_2\) in \(G\) except for \(u_{1}\). Since \(k \geq 2\), there exists at least one vertex \(z_1
\in M_3\) such that \(K_{n}^{z_1} \subseteq V(D_{1})\), and one vertex \(w_1 \in M_{4}\) such that $K_{n}^{w_1} \subseteq V(D_{2})$.

If there exists a vertex \(z^{\prime} \in N_{G}(z_1)\) such that \(K_n^{z^{\prime}} \subseteq V(D_1)\), and there exists a vertex \(w^{\prime} \in N_{G}(w_1)\) such that \(K_n^{w^{\prime}} \subseteq V(D_2)\). By Lemma 2.3, we have  \(|S| \geq n(n - 1)\lambda_2(G)\).
So we suppose that \(K_{n}^{z} \cap V(D_{2}) \neq \emptyset\) for  every vertex \(z \in N_{G}(z_{1})\), or
\(K_{n}^{w} \cap V(D_{1}) \neq \emptyset\) for every vertex \(w \in N_{G}(w_{1})\).
Without loss of generality, assume that for every vertex \(z \in N_{G}(z_{1})\), it holds that \(K_{n}^{z} \cap V(D_{2}) \neq \emptyset\).

If every vertex $z\in N_{G}(z_1)\backslash\{u_1\}$ satisfies \( K_n^{z} \subseteq V(D_2) \), then $|[K_n^{z}, K_n^{z_1}]_{\mathcal{G}} \cap S| = n(n-1)$. Thus,
\begin{align*}
|S| &\geq \left|[K_n^{u_1}, K_n^{u_2}]_{\mathcal{G}} \cap S\right| + \sum_{\substack{u \in N_{G}(u_1) \setminus \{u_2\}}} \left|[K_n^{u_1}, K_n^{u}]_{\mathcal{G}} \cap S\right| \\
&+ \sum_{\substack{u \in N_{G}(u_2) \setminus \{u_1\}}} \left|[K_n^{u_2}, K_n^{u}]_{\mathcal{G}} \cap S\right| + \sum_{\substack{z\in N_{G}(z_1) \setminus \{u_1\}}} |[K_n^{z}, K_n^{z_1}]_{\mathcal{G}} \cap S| \\
&\geq 2(n-2) + (k-1)(n-1) + (k-1)(n-1) + n(n-1)(k-1) \\
&= (n-1)(2k + nk - n) - 2 \\
&\geq(n-1)(2k + 5k - 5) - 2 \\
&> 3k(n-1) - 4.
\end{align*}

If there exists a vertex $z' \in N_{G}(z_1)\backslash\{u_1\}$ such that \( K_n^{z'} \cap V(D_1) \neq \emptyset \) and \( K_n^{z'} \cap V(D_2) \neq \emptyset \), then \( z' \notin N_{G}(u_1) \) and \( z' \notin N_{G}(u_2) \). By  $\lambda(K_2 \times K_n) = n-1$, it follows that $\left| [K_n^{z_1},K_n^z]_{\mathcal{G}} \cap S \right| \geq n-1$ for any vertex $z \in N_{G}(z_1)\backslash\{u_1\}$.  Thus,
\begin{align*}
|S| &\geq \left|[K_n^{u_1}, K_n^{u_2}]_{\mathcal{G}} \cap S\right| + \sum_{\substack{u \in N_{G}(u_1) \setminus \{u_2\}}} \left|[K_n^{u_1}, K_n^{u}]_{\mathcal{G}} \cap S\right| \\
&+ \sum_{\substack{u \in N_{G}(u_2) \setminus \{u_1\}}} \left|[K_n^{u_2}, K_n^{u}]_{\mathcal{G}} \cap S\right| + \sum_{\substack{z \in N_{G}(z_1) \setminus \{u_1\}}} |[K_n^{z}, K_n^{z_1}]_{\mathcal{G}} \cap S| \\
&+ \sum_{\substack{z \in N_{G}(z') \setminus \{z_1\}}} |[K_n^{z}, K_n^{z'}]_{\mathcal{G}} \cap S| \\
&\geq 2(n-2) + (k-1)(n-1) + (k-1)(n-1) + (k-1)(n-1) + (k-1)(n-1) \\
&= (n-1)(4k - 2) - 2 \\
&> 3k(n-1) - 4.
\end{align*}

\noindent{\bf Subcase 2.2.3.}\@ There exist two vertices \( u_4,  u_5 \in N_{G}(u_1) \) such that \( K_n^{u_4} \subseteq V(D_1) \) and \( K_n^{u_5} \subseteq V(D_2) \), or there exist two vertices \( u_4', u_5' \in N_{G}(u_2) \) such that \( K_n^{u_4'} \subseteq V(D_1) \) and \( K_n^{u_5'} \subseteq V(D_2) \).

Without loss of generality, assume that there exist two vertices \(u_{4}, u_{5} \in N_{G}(u_{1})\) such that \(K_{n}^{u_{4}} \subseteq V(D_{1})\) and \(K_{n}^{u_{5}} \subseteq V(D_{2})\). Thus, we have $|[K_n^{u_1}, K_n^{u_4}]_{\mathcal{G}} \cap S| + |[K_n^{u_1}, K_n^{u_5}]_{\mathcal{G}} \cap S| = n(n - 1).$

If there exists a vertex \(u^{\prime} \in N_{G}(u_{4})\) such that \(K_{n}^{u^{\prime}} \subseteq V(D_{1})\), and there exists a vertex \(u^{\prime\prime} \in N_{G}(u_{5})\) such that \(K_{n}^{u^{\prime\prime}} \subseteq V(D_{2})\), then by Lemma 2.3, we have $|S| \geq n(n - 1)\lambda_2(G).$
So we suppose that  \(K_{n}^{u} \cap V(D_{2}) \neq \emptyset\) for every vertex \(u \in N_{G}(u_{4})\),
or \(K_{n}^{u} \cap V(D_{1}) \neq \emptyset\) for every vertex \(u \in N_{G}(u_{5})\). Without loss of generality, assume that for every vertex \(u \in N_{G}(u_{4})\), it holds that \(K_{n}^{u} \cap V(D_{2}) \neq \emptyset\). By  $\lambda(K_2 \times K_n) = n-1$, it follows that $\left| [K_n^{u_4}, K_n^u]_{\mathcal{G}} \cap S \right| \geq n-1$ for any vertex $u \in N_{G}(u_{4})$.

When \( u_4 \notin N_{G}(u_2) \), we have
\begin{align*}
	|S| &\geq \left|[K_n^{u_1}, K_n^{u_2}]_{\mathcal{G}} \cap S\right| +  \left|[K_n^{u_1}, K_n^{u_4}]_{\mathcal{G}} \cap S\right| \\ &+\left|[K_n^{u_1}, K_n^{u_5}]_{\mathcal{G}} \cap S\right|
	+
	\sum_{\substack{u \in N_{G}(u_1) \backslash \{u_2, u_4, u_5\}}} \left|[K_n^{u_1}, K_n^{u}]_{\mathcal{G}} \cap S\right| \\
	&+ \sum_{\substack{u \in N_{G}(u_2) \backslash \{u_1\}}} \left|[K_n^{u_2}, K_n^{u}]_{\mathcal{G}} \cap S\right|
	 +\sum_{\substack{u \in N_{G}(u_4) \backslash \{u_1\}}} \left|[K_n^{u_4}, K_n^{u}]_{\mathcal{G}} \cap S\right| \\
	&\geq 2(n-2) +n(n-1)+ (k-3)(n-1) + (k-1)(n-1) + (k-1)(n-1) \\
	&= (n-1)(3k + n - 3) - 2 \\
	&> 3k(n-1) - 6.
\end{align*}

When \( u_4 \in N_{G}(u_2) \), we have
\begin{align*}
	|S| &\geq \left|[K_n^{u_1}, K_n^{u_2}]_{\mathcal{G}} \cap S\right| +  \left|[K_n^{u_1}, K_n^{u_4}]_{\mathcal{G}} \cap S\right| \\ &+\left|[K_n^{u_1}, K_n^{u_5}]_{\mathcal{G}} \cap S\right|
	+
	\sum_{\substack{u \in N_{G}(u_1) \backslash \{u_2, u_4, u_5\}}} \left|[K_n^{u_1}, K_n^{u}]_{\mathcal{G}} \cap S\right| \\
	&+ \sum_{\substack{u \in N_{G}(u_2) \backslash \{u_1\}}} \left|[K_n^{u_2}, K_n^{u}]_{\mathcal{G}} \cap S\right|+
	 \sum_{\substack{u \in N_{G}(u_4) \backslash \{u_1, u_2\}}} \left|[K_n^{u_4}, K_n^{u}]_{\mathcal{G}} \cap S\right| \\
	&\geq 2(n-2) +n(n-1)+ (k-3)(n-1) + (k-1)(n-1) + (k-2)(n-1) \\
	&= (n-1)(3k + n - 4) - 2 \\
	&> 3k(n-1) - 4.
\end{align*}

The proof is thus complete. $\Box$

If we replace \(K_{n}\) with \(T_{n}\) in Theorem 3.2 and adopt a similar argument, we can obtain the 3-restricted edge-connectivity of \(G \times T_{n}\). Since $\lambda_3(K_2 \times T_n) = 3n-4$ by Lemma 2.7, we assume \(k \geq 2\) in the following theorem.

\begin{theorem}\label{5}
Let \( G \) be a $k$ $(\geq 2)$-regular connected graph with at least four vertices. Then
\[
\lambda_{3}(G \times T_{n}) =
\begin{cases}
\min\{n^{2}\lambda_{2}(G), 3nk - 6\}, & \text{if }g(G) = 3; \\
\min\{n^{2}\lambda_{2}(G), 3nk - 4\}, & \text{if } g(G) \geq 4,
\end{cases}
\]
where \( n \geq 3 \).
\end{theorem}

If \(G\) is maximally restricted edge-connected in Theorem 3.1, then \(\lambda_2(G) =\xi(G)= 2(k - 1)\). This implies \(2n\lambda_2(G) = 4n(k - 1)\). Since \(n\) is odd and \(n \geq 3\), we further derive \(2n\lambda_2(G) > 6k - 4\). Thus \(\lambda_{3}(G \times C_{n}) = \xi_{3}(G \times C_{n})\). Similarly,  if  \(G\) is maximally restricted edge-connected in Theorems 3.2 and 3.3, we also have \(\lambda_{3}(G \times K_{n}) = \xi_{3}(G \times K_{n})\) and \(\lambda_{3}(G \times T_{n}) = \xi_{3}(G \times T_{n})\). The following three corollaries are thus obtained.

\begin{corollary}
Let \(G\) be a $k$ $(\geq 2)$-regular graph \(G\) and $n$ be an odd integer with \(n\geq3\). If \(G\) is maximally restricted edge-connected with \(|V(G)|\geq4\), then \(G \times C_{n}\) is maximally 3-restricted edge-connected.
\end{corollary}

\begin{corollary}
	Let \(G\) be a $k$ $(\geq 2)$-regular graph \(G\) and \(n\geq5\). If \(G\) is maximally restricted edge-connected with \(|V(G)|\geq4\), then \(G \times K_{n}\) is maximally 3-restricted edge-connected.
\end{corollary}

\begin{corollary}
	Let \(G\) be a $k$ $(\geq 2)$-regular graph \(G\) and \(n\geq3\). If \(G\) is maximally restricted edge-connected with \(|V(G)|\geq4\), then \(G \times T_{n}\) is maximally 3-restricted edge-connected.
\end{corollary}

\section{Concluding Remarks}
In this paper, we determine the 3-restricted edge-connectivity for the direct product graphs of a connected regular graph with the cycle, the complete graph and the total graph. Furthermore, we prove that if this  regular graph is maximally restricted edge-connected, then the corresponding direct product graphs are maximally 3-restricted edge-connected. Future research may focus on  the 3-restricted edge-connectivity for the direct product graphs of a connected general graph  with the cycle, the complete graph and the total graph.



\end{sloppypar}
\end{document}